\documentclass{ifacconf}

\usepackage{graphicx}      
\usepackage{natbib}        

\usepackage{amsmath}
\usepackage{enumerate}

\usepackage{xcolor}

\usepackage{todo}

\newcommand{\cC}{\mathcal{C}}
\newcommand{\cH}{\mathcal{H}}
\newcommand{\cL}{\mathcal{L}}
\newcommand{\cO}{\mathcal{O}}
\newcommand{\cP}{\mathcal{P}}
\newcommand{\D}{\mathrm{D}}

\renewcommand{\d}{\mathrm{d}}

\newcommand{\R}{\mathbb{R}}

\newcommand{\der}[2]{\frac{\partial #1}{\partial #2}}

\begin{document}
\pagestyle{headings}

\begin{frontmatter}

\vspace{-1cm}
\emph{This work has been submitted to IFAC for possible publication}

\title{Superconvergence of Galerkin variational integrators
	\thanksref{footnoteinfo}} 

\thanks[footnoteinfo]{MV is supported by the DFG Research Fellowship (VE 1211/1-1)}

\author[First]{Sina Ober-Blöbaum} 
\author[Second]{Mats Vermeeren} 

\address[First]{Department of Mathematics, University of Paderborn, Paderborn, Germany (e-mail: sinaober@math.uni-paderborn.de).}
\address[Second]{School of Mathematics, University of Leeds, Leeds, LS2 9JT, UK (e-mail: m.vermeeren@leeds.ac.uk)}

\begin{abstract}                
We study the order of convergence of Galerkin variational integrators for ordinary differential equations. Galerkin variational integrators approximate a variational (Lagrangian) problem by restricting the space of curves to the set of polynomials of degree at most $s$ and approximating the action integral using a quadrature rule. We show that, if the quadrature rule is sufficiently accurate, the order of the integrators thus obtained is $2s$.
\end{abstract}


\end{frontmatter}

\section{Introduction}

Variational integrators are a class of geometric integration methods, constructed using a discrete version of Hamilton's principle. Variational integrators are symplectic and momentum preserving, provided the discretization exhibits the same symmetries as the continuous system. This leads to more accurate results, especially in long-term simulations of conservative systems (see e.g.~\cite{marsden2001discrete,hairer2006geometric,Reic94}) but also in the simulation and optimization of dissipative systems (\cite{Kane00,Modin2011,JO17,LOH18}) and controlled systems (\cite{ObJuMa10}).

To construct higher order methods, Galerkin variational integrators were considered by \cite{marsden2001discrete} and analyzed by \cite{leok2012general,hall2015spectral,Ca13,ober2015construction} for classical conservative systems. They were further studied by \cite{Wenger2017} for constrained systems and by \cite{COT14} for optimally controlled systems.
The construction of a Galerkin variational integrator consists of two steps:  (1) the approximation of the space of trajectories by a finite-dimensional function space, which we take to be the space of polynomials of degree at most $s$, and (2) the approximation of the action integral by an appropriate quadrature rule.
Obviously, more accurate Galerkin variational integrators are obtained by considering larger function spaces (higher degree of polynomials) and more accurate quadrature rules. We seek to determine the order of a Galerkin variational integrator as a function of the degree of the polynomial approximation and the order of the quadrature rule.

In \cite{hall2015spectral} it is proved that the order of a Galerkin variational integrator based on polynomials of degree at most $s$ and a quadrature rule of order $u$ is at least $\min(s,u)$. Hence, if the quadrature rule is sufficiently accurate, the degree of the polynomials is a lower bound for the order of the integrator. 
Numerical studies in \cite{ober2015construction} indicate that Galerkin variational integrators based on the Lobatto and Gauss quadrature rules are of order $\min(2s,u)$. Hence for a sufficiently accurate quadrature rule, the order of the Galerkin variational integrator seems to be \emph{twice} the degree of polynomials.
For particular classes of Galerkin variational integrators that are equivalent to (modified) symplectic Runge-Kutta methods, a proof of this superconvergence result is provided in \cite{O14}.

In this paper we provide a general proof of the superconvergence of Galerkin variational integrators (Theorem~\ref{thm-superh}). We use techniques similar to \cite{hall2015spectral} and in addition use aspects of the calculus of variations to improve the estimate of the numerical error in the action functional.

\section{Overview of continuous and discrete Lagrangian mechanics}

\subsection{Continuous Lagrangian mechanics}
We consider Lagrange functions $\cL: T Q \rightarrow \R$ on a vector space $Q$. The \emph{action} of a smooth curve $[a,b] \rightarrow Q$ is given by
\[ \mathfrak{S}[q] = \int_a^b \cL(q(t),\dot{q}(t)) \,\d t .\]
A curve $q$ is a \emph{stationary curve} of the action functional $\mathfrak{S}$ if its \emph{Gateaux derivative} (often called \emph{first variation} in this context),
\[ \d \mathfrak{S}[q,\delta q] = \der{}{\alpha}  \mathfrak{S} [q + \alpha \delta q]\Big|_{\alpha = 0} ,\]
vanishes for all smooth curves $\delta q: [a,b] \rightarrow Q$ with $\delta q(a) = \delta q(b) = 0$. Note that $\delta q$ does not denote an infinitesimal, but rather a curve indicating the direction of variation.

The stationarity condition is known as \emph{Hamilton's principle}. It can be expressed as a differential equation using integration by parts:
\begin{align*}
0 = \d \mathfrak{S}[q,\delta q]
&= \int_a^b \der{\cL}{q} \delta q + \der{\cL}{\dot{q}} \delta \dot{q} \,\d t \\
&= \int_a^b \left( \der{\cL}{q} - \frac{\d}{\d t} \der{\cL}{\dot{q}} \right) \delta q \,\d t
+ \left.\left(\der{\cL}{\dot{q}} \delta q \right)\right|_a^b .
\end{align*}
Since $\delta q$ vanishes at the endpoints, the boundary term is zero. The integral is zero for all variations $\delta q$ if and only the \emph{Euler-Lagrange equation}
\[ \der{\cL}{q} - \frac{\d}{\d t} \der{\cL}{\dot{q}} = 0 \]
is satisfied.

Throughout this work we will assume that the Lagrangian is  non-degenerate, i.e.\@ $\det \left( \der{^2 \cL}{\dot{q}^2} \right) \neq 0$. Then the Euler-Lagrange equation is a second order ODE,
\[ \ddot{q} = \left( \der{^2 \cL}{\dot{q}^2} \right)^{-1} \left( \der{\cL}{q} - \der{^2\cL}{\dot{q} \partial q} \dot{q} \right) . \]
In addition, non-degeneracy implies that the Legendre transform $TQ \rightarrow T^*Q: (q,\dot q) \mapsto (q,p) = \left(q, \der{\cL}{\dot{q}} \right)$ is invertible, so the Euler-Lagrange equation is equivalent to the Hamiltonian system
\[ \dot{q} = \der{\cH(q,p)}{p} , \qquad \dot{p} = -\der{\cH(q,p)}{q} \]
with Hamiltonian
\begin{equation}\label{Hamiltonian}
\cH(q,p) = p \dot{q} - \cL(q,\dot{q}) .
\end{equation} 

It is important to note that while the Euler-Lagrange equation is a necessary condition for $q$ to be a minimizer of the action, it is not in general a sufficient condition. In the present work we will make the assumption that solutions to the Euler-Lagrange equations are always minimizers of the action (for sufficiently small intervals of integration). In particular, this is the case for Lagrangians of mechanical type, where $Q = \R^n$ and
\[ \cL = \dot{q}^T M(q) \dot{q} - V(q), \]
where $M$ is a positive definite $n \times n$ matrix \cite[Chapter 5]{gelfand1963calculus}.

\subsection{Variational integrators}
\label{sec-integrators}

Lagrangian systems have a rich structure: their flows consist of symplectic maps (indeed they are equivalent to Hamiltonian systems) and symmetries of the action correspond to conserved quantities by Noether's theorem. When approximating Lagrangian systems numerically, preserving this structure generally leads to improved numerical behaviour. Such an approach of structure-preserving discretization is known as \emph{geometric numerical integration}, see \cite{hairer2006geometric}. In the case of Lagrangian systems, the key to geometric numerical integration is to discretize the action functional instead of discretizing the Euler-Lagrange equation directly. Numerical methods obtained in this way are called \emph{variational integrators}. Below we present some essential facts on variational integrators. For a more detailed discussion we refer to \cite{marsden2001discrete}.

A \emph{discrete Lagrange function} on $Q$ is a differentiable function $L: Q \times Q \times (0,\infty) \rightarrow \R$. The \emph{discrete action} corresponding to a discrete curve $q = (q_0,q_1,\ldots,q_N)$ with step size $h$ is given by 
\[ \mathfrak{S}_d(q) = \sum_{i=1}^N L(q_{i-1}, q_i; h) . \]
We say that $q = (q_0,q_1,\ldots,q_N)$ is \emph{critical} if 
\[ \der{\mathfrak{S}_d(q)}{q_i} = 0 \qquad \text{for } i \in \{1,\ldots,N-1\}, \]
i.e.\@ if the action is invariant with respect to infinitesimal variations of the interior points. This is the case if and only if $q$ satisfies the discrete Euler-Lagrange equation
\begin{equation}\label{dEL}
\D_2 L (q_{i-1}, q_i; h) + \D_1 L(q_i, q_{i+1}; h)  = 0
\end{equation}
for $i \in \{1,\ldots,N-1\}$, where $\D_1$ and $\D_2$ denote the partial derivatives of $L$ with respect to the first and second entry. 

The discrete Euler-Lagrange equation can be interpreted as equality of the two formulas for the discrete momentum,
\[ p_i = \D_2 L(q_{i-1}, q_i; h) \]
and
\[ p_i = -\D_1 L(q_i, q_{i+1}; h) . \]
This gives a natural implementation of the discrete Euler-Lagrange equation as a one-step method
\begin{equation}\label{symplecitc}
\Phi_h: (q_i,p_i) \mapsto (q_{i+1},p_{i+1}) , 
\end{equation}
which is a symplectic integrator for the Hamiltonian system \eqref{Hamiltonian}.

\subsection{Variational error analysis}
\label{sec-varerrana}

Here and in the following, we assume that for sufficiently small $h$ and any pair of boundary values $(q_0,q_1)$ there exists a unique smooth minimizer $q:[0,h] \rightarrow Q$ of the action
\[ \int_0^h \cL(q,\dot{q}) \,\d t, \]
subject to $q(0) = q_0$ and $q(h) = q_1$.
The corresponding minimal value of the action is called the \emph{exact discrete Lagrangian} and denoted by
\[ L_{exact}(q_0,q_1,h) = \min_{\substack{q \in \cC^\infty([0,h],Q) \\ q(0) = q_0, q(h) = q_1}} \int_0^h \cL(q,\dot{q}) \,\d t.
\]
The order of a variational integrator can be determined by comparing its discrete Lagrangian to the exact discrete Lagrangian.
\begin{thm}\label{thm-varerrana}
If for every smooth curve $q$ there holds that
\[ L(q(0),q(h),h) - L_{exact}(q(0),q(h),h) = \cO(h^{\ell+1}), \]
then the variational integrator $\Phi_h$ defined by the discrete Lagrangian $L$, in its symplectic form \eqref{symplecitc}, is of order $\ell$, i.e.\@ it satisfies
\[ \Phi_h(q,p) - \varphi_h(q,p) = \cO(h^{\ell+1}), \]
where $\varphi_h$ is the flow over time $h$ of the continuous Hamiltonian system \eqref{Hamiltonian}.
\end{thm}

This result was first stated in \cite{marsden2001discrete}, where the proof contained a flaw which was later fixed by \cite{patrick2009error}. The difficulty of the proof lies in a singularity in the discrete Legendre transform when $h \rightarrow 0$. In a forthcoming work, we will present a new approach to this result from the perspective of \emph{modified Lagrangians} (\cite{vermeeren2017modified}).

\subsection{Galerkin variational integrators}
\label{galerkin}

An effective method to construct higher order variational integrators is to use a \emph{Galerkin discretization}. To construct a Galerkin integrator, the space of smooth curves on the time interval of one step, $\cC^\infty([0,h],Q)$, is replaced by a finite dimensional space of polynomials
\begin{align*}
\cP^s([0,h],Q) = \{ q \in \cC^\infty([0,h],Q) \,|\, &q \text{ a polynomial} \\ &\text{of degree at most } s \} .
\end{align*}
We fix $s+1$ control points $h d_0 < h d_1 < \ldots < h d_s$, where $d_0 = 0$ and $d_s = 1$. If for each of these control points a value $q(h d_j) = q_j$ is prescribed, then the polynomial $q \in \cP^s([0,h],Q)$ is uniquely determined. We denote by $\hat{q}(\,\cdot\, ; q_0,\ldots,q_s,h)$ the polynomial thus obtained. 

Given a continuous Lagrangian $\cL$, we define
\begin{align*}
&\cL_\mathrm{p}(t ; q_0,\ldots,q_s,h) \\
&\qquad = \cL \left( \hat{q}(t ; q_0,\ldots,q_s,h), \dot{\hat{q}}(t ; q_0,\ldots,q_s,h) \right)
\end{align*}
where the subscript p reminds us that this is the Lagrangian evaluated on a polynomial. We would like to consider the discrete action
\[ \int_0^h \cL_\mathrm{p}(t ; q_0,\ldots,q_s,h) \,\d t. \]
However, to evaluate this integral numerically, we need a quadrature rule. We fix quadrature points $c_i \in [0,1]$ and weights $b_i \in \R$, with $\sum_i b_i = 1$. We denote by $u$ the order of the corresponding quadrature rule. Then for any smooth function $f$ there holds
\[ \int_0^h f(t) \,\d t - h \sum_i b_i f(h c_i) = \cO(h^{u+1}). \]
We define the discrete Lagrangian as
\[ L(q_0,q_s;h) = \min_{q_1,\ldots,q_{s-1} \in Q}
\left( h \sum_i b_i \cL_\mathrm{p}(h c_i ; q_0,\ldots,q_s,h) \right) ,
\]
Here we assume that there exists a unique minimizer. This is the case in particular if the Lagrangian is of mechanical type and the quadrature rule is sufficiently accurate \cite[Theorem 3.5]{hall2015spectral}.

The discrete Lagrangian can also be written as
\[
L(q_0,q_s;h) 
= \min_{\substack{q \in \cP^s([0,h],Q), \\ q(0) = q_0, q(h) = q_s}}  \mathfrak{S}_\mathrm{i}[q] .
\]
where $\mathfrak{S}_\mathrm{i}$ denotes the \emph{internal action},
\begin{equation}\label{quadrature-action}
\mathfrak{S}_\mathrm{i}[q] = h \sum_i b_i \cL( q(h c_i), \dot{q} (h c_i) ) .
\end{equation}

For more details on the construction of Galerkin variational integrators, see for example \cite{marsden2001discrete,leok2012general,ober2015construction}.

\section{A few technicalities}

Before we can prove our main result on the superconvergence of Galerkin variational integrators, we need some error estimates for polynomial interpolation.

\begin{lem}\label{lemma-poly}
	Let $q$ be smooth curve and $\hat{q}$ a family of polynomials of degree $s$, parametrized by $h$, which equals $q$ at the control points $0 = h d_0 < h d_1 < \ldots < h d_s = h$.
	Then for any $k \leq s$ there holds
	\[ \big\| q^{(k)} - \hat{q}^{(k)} \big\|_\infty = \cO(h^{s+1-k}) , \]
	where $\|\cdot\|_\infty$ denotes the maximum norm on $[0,h]$.
\end{lem}
\begin{pf}
	Since $q- \hat{q}$ has at least $s+1$ zeros in the interval $[0,h]$, we know by the mean value theorem that $\dot{q} - \dot{\hat{q}}$ has at least $s$ zeros, and recursively we find that $q^{(k)}- \hat{q}^{(k)}$ has at least $s+1-k$ zeroes. In particular, $q^{(s)}- \hat{q}^{(s)}$ has at least one zero $t_0 \in [0,h]$. Since $\hat{q}^{(s+1)} = 0$ identically, it follows that:
	\begin{align*}
	\big\| q^{(s)} - \hat{q}^{(s)} \big\|_\infty
	&\leq \big\| q^{(s+1)}\big \|_\infty \max_{t \in [0,h]} (t - t_0) \\
	&= \cO(h) . 
	\end{align*}
	Combining this with the fact that $q^{(s-1)} - \hat{q}^{(s-1)}$ has a zero $t_1$ in $[0,h]$, we find
	\begin{align*}
	\big\| q^{(s-1)} - \hat{q}^{(s-1)} \big\|_\infty 
	&\leq \big\| q^{(s)} - \hat{q}^{(s)} \big\|_\infty \max_{t \in [0,h]} (t - t_1) \\
	&= \cO(h^2) . 
	\end{align*}
	Repeating this argument recursively, we obtain
	\[ \big\| q^{(k)} - \hat{q}^{(k)} \big\|_\infty = \cO(h^{s+1-k}) . \qed \]
\end{pf}

The following proposition contains some simple inequalities that will be useful below. 

\begin{prop}\label{prop-inequalities}
	For any differentiable curve $\delta q: [0,h] \rightarrow Q$ with $\delta q(0) = \delta q(h) = 0$ there holds
	\begin{enumerate}[$(a)$]
		\item $\displaystyle \| \delta q \|_\infty \leq \frac{1}{2} \| \delta \dot{q} \|_1$ ,
		\item $\displaystyle \| \delta q \|_1 \leq \frac{h}{2} \| \delta \dot{q} \|_1$ ,
	\end{enumerate}
	where $\| \cdot \|_p$ denotes the $L^p$-norm on $[0,h]$.
	Furthermore, for any differentiable curve $q: [0,h] \rightarrow Q$ there holds
	\begin{enumerate}[$(c)$]
		\item $\| q \|_1 \leq \sqrt{h} \| q \|_2$ .
	\end{enumerate}
\end{prop}
\begin{pf}
	\begin{enumerate}[$(a)$]
		\item Let $|\delta q|$ reach its maximum in $(0,h)$ at $t_{max}$. We have
		\begin{align*} \| \delta \dot{q} \|_1
		&= \int_0^{t_{max}} |\delta \dot{q}| \,\d t + \int_{t_{max}}^h |\delta \dot{q}| \,\d t \\
		&\geq |\delta q(t_{max}) - \delta q(0)| + |\delta q(h) - \delta q(t_{max})| \\
		&= 2 \| \delta q \|_\infty . 
		\end{align*}
		\item We have
		\[ \| \delta q \|_1 = \int_0^h |\delta q(t)| \,\d t \leq \int_0^h \|\delta q\|_\infty \,\d t = h \|\delta q\|_\infty , \]
		so the claim follows from inequality $(a)$.
		\item This is a special case of Hölder's inequality, 
		\[ \|f g\|_1 \leq \|f\|_\alpha \|g\|_\beta, \]
		with $f = 1$, $g = q$ and $\alpha = \beta = 2$. \qed
	\end{enumerate}
\end{pf}

\section{Superconvergence}
\label{sec-super}

We now come to our main result.

\begin{thm}
	\label{thm-superh}
	Let $L$ be a Galerkin discretization of a Lagrangian $\cL$, based on polynomials of degree $s$ and a quadrature rule of degree $u$. Assume that all discrete and continuous critical curves minimize their respective actions. 
	Then the corresponding symplectic integrator \eqref{symplecitc} is of order $\min(2 s,u)$.
\end{thm}

\begin{pf}
	By Theorem \ref{thm-varerrana}, it suffices to show that
	\[ L_\mathrm{exact}(q(0),q(h);h) - L(q(0),q(h);h) = \cO \left( h^{\min(2 s,u) + 1} \right) \]
	for every smooth curve $q$.
	
	Let $q_\mathrm{EL}$ denote the unique minimizer of the continuous action with $q_\mathrm{EL}(0) = q(0)$ and $q_\mathrm{EL}(h) = q(h)$. The subscript reminds us that $q_\mathrm{EL}$ satisfies the continuous Euler-Lagrange equation. Let $\hat{q} \in \cP^s([0,h],Q)$ be the polynomial that agrees with $q_\mathrm{EL}$ at the control points $0 = h d_0 < h d_1 < \ldots < h d_s = h$ and $\tilde{q} \in \cP^s([0,h],Q)$ the polynomial that minimizes the internal action
	\[ \mathfrak{S}_\mathrm{i}[q] = h \sum_i b_i \cL(q(h c_i), \dot{q}(h c_i) ) \]
	in $\cP^s([0,h],Q)$. Since
	\[ L_\mathrm{exact}(q_0,q_h;h) = \int_0^h \cL(q_\mathrm{EL},\dot{q}_\mathrm{EL}) \,\d t  \]
	we have to show that
	\begin{equation}\label{toshow}
	\begin{split}
	\int_0^h \cL(q_\mathrm{EL},\dot{q}_\mathrm{EL}) \,\d t - h \sum_i b_i \cL(\tilde{q}(h c_i), \dot{\tilde{q}}(h c_i) ) \\
	= \cO \left(h^{\min(2s,u) + 1}\right) .
	\end{split}
	\end{equation}
	We expand this difference as
	\begin{equation}
	\label{difference-expansion}
	\begin{split}
	& \left( \int_0^h \cL(q_\mathrm{EL}, \dot{q}_\mathrm{EL}) \,\d t - \int_0^h \cL(\hat q, \dot{\hat q}) \,\d t \right)  \\
	& + \left( \int_0^h \cL(\hat q, \dot{\hat q}) \,\d t - h \sum_i b_i \cL(\tilde{q}(h c_i), \dot{\tilde{q}}(h c_i) ) \right) .
	\end{split}
	\end{equation}
	
	We start with the first term of \eqref{difference-expansion}. From Lemma \ref{lemma-poly} we know that $q - \hat{q} = \cO(h^{s+1})$ and $\dot{q} - \dot{\hat{q}} = \cO(h^s)$, hence
	\begin{align*}
	& \int_0^h \cL(q_\mathrm{EL}, \dot{q}_\mathrm{EL}) \,\d t - \int_0^h \cL(\hat q, \dot{\hat q}) \,\d t \\
	&= \int_0^h \bigg( \der{\cL(q_\mathrm{EL}, \dot{q}_\mathrm{EL})}{q}(q_\mathrm{EL} - \hat q) +  \der{\cL(q_\mathrm{EL}, \dot{q}_\mathrm{EL})}{\dot{q} } (\dot{q}_\mathrm{EL} - \dot{\hat q})\\
	&\hspace{6.5cm} + \cO\!\left(h^{2s}\right) \bigg) \d t \\
	&= \int_0^h \bigg( \left(\der{\cL(q_\mathrm{EL}, \dot{q}_\mathrm{EL})}{q} - \frac{\d}{\d t} \der{\cL(q_\mathrm{EL}, \dot{q}_\mathrm{EL})}{\dot{q} } \right)(q_\mathrm{EL} - \hat q) \\
	&\hspace{6.5cm} + \cO\!\left(h^{2s}\right) \bigg) \d t \\
	&\qquad + \left( \der{\cL(q_\mathrm{EL}, \dot{q}_\mathrm{EL})}{\dot{q} } (q_\mathrm{EL} - \hat q) \middle) \right|_0^h .
	\end{align*}
	The boundary term vanishes because $\hat q(0) = q_\mathrm{EL}(0)$ and $\hat q(h) = q_\mathrm{EL}(h)$. Furthermore, $q_\mathrm{EL}$ solves the Euler-Lagrange equation, so we find
	\begin{equation}\label{difference1}
	\begin{split}
	\int_0^h \cL(q_\mathrm{EL}, \dot{q}_\mathrm{EL}) \,\d t - \int_0^h \cL(\hat q, \dot{\hat q}) \,\d t 
	&=  \int_0^h \cO\!\left( h^{2s} \right) \,\d t \\
	&= \cO\!\left( h^{2s+1} \right) . 
	\end{split}
	\end{equation}
	
	To bound the second term of \eqref{difference-expansion} we follow the arguments of \cite[Theorem 3.3]{hall2015spectral}. Since $\tilde{q}$ is the minimizing element of $\cP^s([0,h],Q)$, we have
	\begin{align*}
	h \sum_i b_i \cL(\tilde{q}(h c_i), \dot{\tilde{q}}(h c_i) ) 
	&\leq h \sum_i b_i \cL(\hat{q}(h c_i), \dot{\hat{q}}(h c_i) ) \\
	&\leq \int_0^h \cL(\hat q, \dot{\hat q}) \,\d t + \cO(h^{u+1})
	\end{align*}
	On the other hand, since $q_\mathrm{EL}$ minimizes the continuous action, there holds
	\begin{align*}
	h \sum_i & b_i \cL(\tilde{q}(h c_i), \dot{\tilde{q}}(h c_i) ) \\
	&= \int_0^h \cL(\tilde q, \dot{\tilde q}) \,\d t + \cO(h^{u+1}) \\
	&\geq \int_0^h \cL(q_\mathrm{EL}, \dot{q}_\mathrm{EL}) \,\d t + \cO(h^{u+1}) \\
	&= \int_0^h \cL(\hat q, \dot{\hat q}) \,\d t + \cO(h^{2s+1} + h^{u+1}) ,
	\end{align*}
	where the last line follows from \eqref{difference1}. Combining both inequalities we find
	\begin{equation}\label{difference2}
	\int_0^h \cL(\hat q, \dot{\hat q}) \,\d t - h \sum_i b_i \cL(\tilde{q}(h c_i), \dot{\tilde{q}}(h c_i) ) = \cO(h^{\min(2s,u) + 1}) .
	\end{equation}
	Equations \eqref{difference1} and \eqref{difference2} together imply the desired result \eqref{toshow}.
	\qed
\end{pf}

\section{Convergence of the Galerkin curves}

Theorem \ref{thm-superh} states that for a sufficiently accurate quadrature rule, the one-step method obtained by Galerkin discretization has order $2s$, twice the degree of polynomials used. If we compare the polynomial approximations to the exact solution at arbitrary times (away from the mesh points), we find an error of order $s$, the same as the degree of polynomials. 
This halving of the order was also observed in \cite[Section 3.4]{hall2015spectral}.
Below we prove this claim under a coercivity assumption \eqref{coercivity}. This assumption is satisfied in particular for mechanical Lagrangians, as shown for the case of a constant mass matrix in \cite{hall2015spectral}. We will provide a more general proof in a forthcoming publication.

\begin{thm}\label{thm-galerkin}
	Assume that there exists a $C>0$ such that for every continuous critical curve $q_{EL}$ of the action $\mathfrak{S}$ and for any variation $\delta q$, vanishing at the endpoints, there holds
	\begin{equation}\label{coercivity}
	\mathfrak{S} [q_{EL} + \delta q] - \mathfrak{S}[q_{EL}] \geq C \| \delta \dot{q} \|_2^2 .
	\end{equation}
	Then, for sufficiently small $h$, the polynomial $\tilde{q}$ of degree $s$ minimizing the discrete action satisfies
	\[ \| \tilde{q} - q_{EL} \|_\infty = \cO \left( h^{\min(s,\frac{u}{2}) +  1} \right) \]
	and
	\[ \| \tilde{q} - q_{EL} \|_1 = \cO \left( h^{\min(s,\frac{u}{2}) + 2} \right) \]
	where $\| \cdot \|_p$ denotes the $L^p$-norm on $[0,h]$.
\end{thm}
\begin{pf}
	From the proof of Theorem \ref{thm-superh} we know that
	\begin{align*}
	& \mathfrak{S}[\tilde{q}] - \mathfrak{S}[q_{EL}]\\
	&= \int_0^h \cL(\tilde{q},\dot{\tilde{q}}) - \int_0^h \cL(q_{EL},\dot{q}_{EL}) \,\d t \\
	&= L(q(0),q(h);h) - L_\mathrm{exact}(q(0),q(h);h) + \cO(h^{u+1}) \\
	&= \cO( h^{\min(2s,u) + 1}) .
	\end{align*}
	From \eqref{coercivity}, with $\delta q = \tilde{q} - q_{EL}$ it now follows that
	\[ C \| \dot{\tilde{q}} - \dot{q}_{EL} \|_2^2 = \cO \left( h^{\min(2s,u) + 1} \right)\]
	and hence
	\[ \| \dot{\tilde{q}} - \dot{q}_{EL} \|_2 = \cO \left( h^{\min(s,\frac{u}{2}) + \frac{1}{2}} \right) . \]
By Hölder's inequality (Proposition \ref{prop-inequalities}$(c)$) it follows that 
\[ \| \dot{\tilde{q}} - \dot{q}_{EL} \|_1 = \cO \left( h^{\min(s,\frac{u}{2}) + 1} \right) . \]
Finally, since $q_{EL}(0) = \tilde{q}(0)$ and $q_{EL}(h) = \tilde{q}(h)$, it follows from Proposition \ref{prop-inequalities}$(a)$ that
\[ \| \tilde{q} - q_{EL} \|_\infty = \cO \left( h^{\min(s,\frac{u}{2}) +  1} \right) \]
and from Proposition \ref{prop-inequalities}$(b)$ that
\[ \| \tilde{q} - q_{EL} \|_1 = \cO \left( h^{\min(s,\frac{u}{2}) + 2} \right) . \qed \]
\end{pf}

\section{Possible extension to forced systems}

Lagrangian systems with external forces are an important extension of the theory, especially towards the study of optimal control problems. Variational integrators for systems with external forces were presented in \cite{marsden2001discrete}, along with a brief argument suggesting variational error analysis is possible in this case too. This is studied in detail in a recent preprint by \cite{fernandez2021error}. In the forced case order estimates can be obtained by comparing the discrete Lagrangian to the exact discrete Lagrangian and the discrete forces to the exact discrete forces. Alternatively, a forced system can be embedded into a Lagrangian system without external forces of higher dimension. This approach to variational error analysis of forced systems was taken by \cite{dediego2018variational}.

\subsection{Forced Galerkin integrators}

Forced Lagrangian systems are defined by the Lagrange-d'Alembert principle
\[ \delta \int_a^b \cL(q,\dot{q}) \,\d t + \int_a^b f(q,\dot{q}) \delta q \,\d t = 0. \]
The corresponding (forced) Euler-Lagrange equation is
\begin{equation}\label{fEL}
\der{\cL(q,\dot{q})}{q} - \frac{\d}{\d t}\der{\cL(q,\dot{q})}{\dot{q}} 
+ f(q,\dot{q}) = 0 .
\end{equation}
The discrete Lagrange-d'Alembert principle requires a discrete Lagrangian $L(q_{k-1},q_k;h)$ and discrete forces $F^\pm(q_{k-1},q_k;h)$. It reads
\begin{align*}
\delta \sum_k L(q_k,q_{k+1}) + \sum_k \left( F^-(q_k,q_{k+1}) + F^+(q_{k-1},q_k) \right) \delta q_k \\
= 0
\end{align*}
and yields the equations
\begin{align*}
&\D_1 L(q_k,q_{k+1}) + \D_2 L(q_{k-1},q_k) \\
&\quad + F^-(q_k,q_{k+1}) + F^+(q_{k-1},q_k) = 0.
\end{align*}

Galerkin integrators for forced systems are constructed as follows (see \cite{COT14}). As in Section \ref{galerkin} we denote by $\hat{q}(\,\cdot\, ; q_0,\ldots,q_s,h)$ the polynomial of degree at most $s$ defined by its values $q_0,\ldots,q_s$ at control points $h d_0 < h d_1 < \ldots < h d_s$, where $d_0 = 0$ and $d_s = 1$. We define
\begin{align*}
&\cL_\mathrm{p}(t ; q_0,\ldots,q_s,h) \\
&\qquad = \cL \left( \hat{q}(t ; q_0,\ldots,q_s,h), \dot{\hat{q}}(t ; q_0,\ldots,q_s,h) \right)
\end{align*}
and
\begin{align*}
&f_\mathrm{p}(t ; q_0,\ldots,q_s,h) \\
&\qquad = f \left( \hat{q}(t ; q_0,\ldots,q_s,h), \dot{\hat{q}}(t ; q_0,\ldots,q_s,h) \right) .
\end{align*}
Consider quadrature points $c_i \in [0,1]$ and weights $b_i \in \R$. Given $q_0$ and $q_s$ we impose
\begin{align*}
&\delta \sum_i b_i \cL_\mathrm{p}(h c_i ; q_0,\ldots,q_s,h) \\
&\qquad  + \sum_i b_i f_\mathrm{p}(h c_i ; q_0,\ldots,q_s,h) \delta \hat{q}(h c_i; q_0,\ldots,q_s,h) = 0,
\end{align*}
where $\delta$ stands for arbitrary variations of the interior control values $q_1,\ldots,q_{s-1}$. This gives us $s-1$ equations 
\begin{equation}\label{forcedint}
\begin{split}
&\sum_i b_i \der{\cL_\mathrm{p}(h c_i ; q_0,\ldots,q_s,h)}{q_k} \\
&\quad + \sum_i b_i f_\mathrm{p}(h c_i ; q_0,\ldots,q_s,h) \der{\hat{q}(h c_i ; q_0,\ldots,q_s,h)}{q_k} = 0,
\end{split}
\end{equation}
where $k = 1,\ldots s-1$. Assuming these $s-1$ equations uniquely determine $q_1, \ldots, q_{s-1}$ as functions of $q_0$ and $q_s$, we can define
\begin{equation}\label{dLdA1}
 L(q_0,q_s;h) = \sum_i hb_i \cL_\mathrm{p}(h c_i ; q_0,\ldots,q_s,h) .
 \end{equation}
and
\begin{equation}
\begin{split}
&F^-(q_0,q_s;h) \\
&= h\sum_i b_i f_\mathrm{p}(h c_i ; q_0,\ldots,q_s,h) \der{\hat{q}(h c_i ; q_0,\ldots,q_s,h)}{q_0}, \\
&F^+(q_0,q_s;h) \\
&= h\sum_i b_i f_\mathrm{p}(h c_i ; q_0,\ldots,q_s,h) \der{\hat{q}(h c_i ; q_0,\ldots,q_s,h)}{q_s}.
\end{split} \label{dLdA3}
\end{equation}
Equations \eqref{dLdA1}--\eqref{dLdA3} are the discrete Lagrangian and discrete forces defining the Galerkin integrator with polynomials of degree $s$ and quadrature rule given by $(b_i,c_i)$.

\subsection{Obstruction to proving superconvergence}

In the presence of external forces, the exact discrete Lagrangian depends not just on the Lagrangian, but also on the external forces. It is obtained by evaluating the action over the interval $[0,h]$ on the solution $q_\mathrm{fEL}$ of the \emph{forced} Euler-Lagrange equation \eqref{fEL}. Hence the estimate for the first term in \eqref{difference-expansion} becomes
\begin{align}
& \int_0^h \cL(q_\mathrm{fEL}, \dot{q}_\mathrm{fEL}) \,\d t - \int_0^h \cL(\hat q, \dot{\hat q}) \,\d t \notag\\
&= \int_0^h \bigg( \left(\der{\cL(q_\mathrm{fEL}, \dot{q}_\mathrm{fEL})}{q} - \frac{\d}{\d t} \der{\cL(q_\mathrm{fEL}, \dot{q}_\mathrm{fEL})}{\dot{q} } \right)(q_\mathrm{fEL} - \hat q) \notag\\
&\hspace{6.5cm} + \cO\!\left(h^{2s}\right) \bigg) \d t \notag\\
&= - \int_0^h f(q_\mathrm{fEL},\dot{q}_\mathrm{fEL}) (q_\mathrm{fEL} - \hat q) \,\d t + \cO\!\left(h^{2s+1}\right) .
\label{forced}
\end{align}
For generic forces $f$ we have $f(q_\mathrm{fEL},\dot{q}_\mathrm{fEL}) = \cO(1)$ and $q_\mathrm{fEL} - \hat q = \cO(h^s)$, so we can only estimate \eqref{forced} by $\cO(h^{s+1})$.  We expect a similar estimate for the difference between the exact discrete forces
\begin{align*}
F^-_{\mathrm{exact}}(q_0,q_s;h) = \int_0^h f(q_\mathrm{fEL}, \dot{q}_\mathrm{fEL}) \der{q_\mathrm{fEL}}{q_0}, \\
F^+_{\mathrm{exact}}(q_0,q_s;h) = \int_0^h f(q_\mathrm{fEL}, \dot{q}_\mathrm{fEL}) \der{q_\mathrm{fEL}}{q_s},
\end{align*}
and their numerical approximations. This means that we cannot prove superconvergence using the forced analogue of variational error analysis.

\subsection{Possible workaround}
\label{sec-forced-workaround}

In the previous subsection we observed that our proof of superconvergence fails in the presence of external forces. This is because our proof requires that the dynamics are given by Hamilton's principle rather than the Lagrange-d'Alembert principle. Still, we expect a superconvergence result to hold for forced systems too. This expectation is based on well-understood low-order methods (e.g.\@ the midpoint rule and Störmer-Verlet method are Galerkin integrators based on linear polynomials, but they are second order methods) as well as preliminary numerical observations for higher-order methods.

A potential way to remedy our proof is the observation that forced systems can also be described by Hamilton's principle if we double the dimension and introduce a variable $Q$, which in the end we will require to be a copy of $q$ (see \cite{galley2013classical}). In particular, we consider the extended Lagrangian
\begin{equation}\label{extendedL}
\begin{split}
\cL^\mathrm{f}(q,Q,\dot{q},\dot{Q}) &= \cL(Q,\dot{Q}) - \cL(q,\dot{q}) \\
&\quad+ \tfrac{1}{2} ( f(Q,\dot{Q}) + f(q,\dot{q}) ) (Q - q).
\end{split}
\end{equation}
Taking variations with respect to $Q$ we find the Euler-Lagrange equation
\begin{equation}\label{extendedEL}
\begin{split}
&\der{\cL(Q,\dot{Q})}{Q} - \frac{\d}{\d t}\der{\cL(Q,\dot{Q})}{\dot{Q}} 
+ \frac{ f(q,\dot{q}) + f(Q,\dot{Q})}{2} \\
&\quad + \der{f(Q,\dot{Q})}{Q} \frac{Q-q}{2} - \frac{\d}{\d t} \left( \der{f(Q,\dot{Q})}{\dot{Q}} \frac{Q-q}{2} \right)
= 0 .
\end{split}
\end{equation}
When we impose $Q = q$ this equation reduces to the familiar forced Euler-Lagrange equation \eqref{fEL}. The same conclusion holds for variations with respect to $q$. As pointed out by \cite{dediego2018variational}, this observation can be used to apply variational error analysis to forced systems. In our present context, we need to show that the forced Galerkin integrator \eqref{forcedint} is equivalent to a Galerkin integrator for the extended system \eqref{extendedEL}.

Consider the Galerkin integrator for \eqref{extendedEL} defined by the Lagrangian
\begin{align*}
&\cL_\mathrm{p}^\mathrm{f}(t ; (q_0, Q_0),\ldots,(q_s, Q_s),h) \\
&\quad = \cL^\mathrm{f} \Big( \hat{q}(t ; q_0,\ldots,q_s,h), \hat{q}(t ; Q_0,\ldots,Q_s,h),\\
&\qquad\qquad\quad \dot{\hat{q}}(t ; q_0,\ldots,q_s,h) , \dot{\hat{q}}(t ; Q_0,\ldots,Q_s,h) \Big) \\
&\quad = \cL_\mathrm{p}(t ; Q_0,\ldots,Q_s,h) - \cL_\mathrm{p}( t ; q_0,\ldots,q_s,h) \\
&\qquad + \tfrac{1}{2} \left( f_\mathrm{p}(t ; Q_0,\ldots,Q_s,h) + f_\mathrm{p}( t ; q_0,\ldots,q_s,h) \right) \\
&\qquad\qquad \cdot \left( \hat{q}(t ; Q_0,\ldots,Q_s,h) - \hat{q}(t ; q_0,\ldots,q_s,h) \right) 
\end{align*}
and the quadrature rule with points $c_i$ and weights $b_i$. By Theorem \ref{thm-superh} this integrator is of order $2s$ if the quadrature rule is sufficiently accurate and if all critical curves minimize the action. Varying $Q_k$ (or $q_k$) for some k with $0 < k < s$, and then imposing $Q_\ell = q_\ell$ for all $0 < \ell < s$, leads to the internal equations \eqref{forcedint}. As before, we assume that these uniquely determine $q_1 = Q_1, \ldots, q_{s-1}=Q_{s-1}$ as functions of $q_0,Q_0,q_s,Q_s$, allowing us to define the extended discrete Lagrangian
\begin{align*}
& L^\mathrm{f}(q_0,Q_0,q_s,Q_s,h) \\
&\qquad = \sum_i hb_i \cL_\mathrm{p}^\mathrm{f}(h c_i ; (q_0,Q_0),\ldots,(q_s,Q_s),h) .
\end{align*}
Its discrete Euler-Lagrange equations, evaluated on $Q_0 = q_0$, $Q_s = q_s$, are equivalent to the discrete Lagrange-d'Alembert principle for \eqref{dLdA1}--\eqref{dLdA3}. Hence the forced Galerkin integrator defined by \eqref{dLdA1}--\eqref{dLdA3} is of order $2s$.

The attentive reader may have noticed a problem with the argument above. To apply Theorem \ref{thm-superh} to Galerkin integrator for the extended system, we need te property that critical curves of the action are minimizers, but this does not hold for Lagrangians of the form \eqref{extendedL}. In the proof of Theorem \ref{thm-superh} we used this assumption to show that the minimizing polynomial $\tilde{q}$ is close to the polynomial interpolating the continuous solution $\hat{q}$ and hence to estimate the second term in \eqref{difference-expansion}. However, it is plausible that even without this assumption the difference $\hat{q}-\tilde{q}$ will be small in a generic case. We currently do not have a precise statement of this claim, so this proof is left to be finished in future work.

\section{Conclusion and outlook}

Following the approach of \cite{hall2015spectral}, but with stronger error bounds obtained from the calculus of variations, we have shown that Galerkin variational integrators exhibit superconvergence: given a suitably accurate quadrature rule, the order of such an integrator is \emph{twice} the degree of polynomials used to construct it. In our presentation here we have relied heavily on established results concerning variational error analysis and kept technical details to a minimum. Since some of these details are worthy of attention, we plan to continue this topic in a forthcoming paper, where we will 
\begin{itemize}
	\item Present a new proof of Theorem \ref{thm-varerrana} based on modified Lagrangians. 
	\item Show that mechanical Lagrangians (with a possibly position-dependent mass matrix) and their discretizations satisfy a coercivity condition as assumed in Theorem \ref{thm-galerkin}. From this condition it also follows that critical curves are minimizers.
	\item Provide numerical experiments to illustrate Theorems \ref{thm-superh} and \ref{thm-galerkin}, beyond what is already available in the work of \cite{ober2015construction}.
\end{itemize}
An important additional topic for future work is to turn the arguments sketched in \ref{sec-forced-workaround} into a rigorous proof of superconvergence in the presence of external forces.

\bibliography{gni-wo-doi}

\begin{thebibliography}{21}
\providecommand{\natexlab}[1]{#1}
\providecommand{\url}[1]{\texttt{#1}}
\providecommand{\urlprefix}{URL }
\expandafter\ifx\csname urlstyle\endcsname\relax
  \providecommand{\doi}[1]{doi:\discretionary{}{}{}#1}\else
  \providecommand{\doi}{doi:\discretionary{}{}{}\begingroup
  \urlstyle{rm}\Url}\fi

\bibitem[{Campos(2014)}]{Ca13}
Campos, C.M. (2014).
\newblock High order variational integrators: A polynomial approach.
\newblock In F.~Casas and V.~Martínez (eds.), \emph{Advances in Differential
  Equations and Applications}, volume~4 of \emph{SEMA SIMAI Springer Series},
  249--258. Springer International Publishing.

\bibitem[{Campos et~al.(2015)Campos, Ober-Bl\"obaum, and Tr\'{e}lat}]{COT14}
Campos, C.M., Ober-Bl\"obaum, S., and Tr\'{e}lat, E. (2015).
\newblock High order variational integrators in the optimal control of
  mechanical systems.
\newblock \emph{Discrete and Continuous Dynamical Systems}, 35(9), 4193--4223.

\bibitem[{De~Diego and de~Almagro(2018)}]{dediego2018variational}
De~Diego, D.M. and de~Almagro, R.S.M. (2018).
\newblock Variational order for forced {L}agrangian systems.
\newblock \emph{Nonlinearity}, 31(8), 3814.

\bibitem[{Fernández et~al.(2021)Fernández, Zurita, and
  Grillo}]{fernandez2021error}
Fernández, J., Zurita, S.G., and Grillo, S. (2021).
\newblock Error analysis of forced discrete mechanical systems.
\newblock \emph{arXiv:2103.11060}.

\bibitem[{Galley(2013)}]{galley2013classical}
Galley, C.R. (2013).
\newblock Classical mechanics of nonconservative systems.
\newblock \emph{Physical review letters}, 110(17), 174301.

\bibitem[{Gelfand and Fomin(1963)}]{gelfand1963calculus}
Gelfand, I.M. and Fomin, S.V. (1963).
\newblock \emph{Calculus of Variations}.
\newblock Prentice-Hall.

\bibitem[{Hairer et~al.(2006)Hairer, Lubich, and Wanner}]{hairer2006geometric}
Hairer, E., Lubich, C., and Wanner, G. (2006).
\newblock \emph{Geometric Numerical Integration: Structure-Preserving
  Algorithms for Ordinary Differential Equations}.
\newblock Springer, 2nd edition.

\bibitem[{Hall and Leok(2015)}]{hall2015spectral}
Hall, J. and Leok, M. (2015).
\newblock Spectral variational integrators.
\newblock \emph{Numerische Mathematik}, 130(4), 681--740.

\bibitem[{Jim\'{e}nez and Ober-Bl\"obaum(2018)}]{JO17}
Jim\'{e}nez, F. and Ober-Bl\"obaum, S. (2018).
\newblock A fractional variational approach for modelling dissipative
  mechanical systems: Continuous and discrete settings.
\newblock \emph{IFAC-PapersOnLine}, 51(3), 50 -- 55.
\newblock 6th IFAC Workshop on Lagrangian and Hamiltonian Methods for Nonlinear
  Control LHMNC 2018.

\bibitem[{Kane et~al.(2000)Kane, Marsden, Ortiz, and West}]{Kane00}
Kane, C., Marsden, J.E., Ortiz, M., and West, M. (2000).
\newblock Variational integrators and the {N}ewmark algorithm for conservative
  and dissipative mechanical systems.
\newblock \emph{International Journal for Numerical Methods in Engineering},
  49(10), 1295--1325.

\bibitem[{Leok and Shingel(2012)}]{leok2012general}
Leok, M. and Shingel, T. (2012).
\newblock General techniques for constructing variational integrators.
\newblock \emph{Frontiers of Mathematics in China}, 7(2), 273--303.

\bibitem[{{Limebeer} et~al.(2020){Limebeer}, {Ober-Bloebaum}, and {Haddad
  Farshi}}]{LOH18}
{Limebeer}, D., {Ober-Bloebaum}, S., and {Haddad Farshi}, F. (2020).
\newblock Variational integrators for dissipative systems.
\newblock \emph{IEEE Transactions on Automatic Control}, 65(4), 1381--1396.

\bibitem[{Marsden and West(2001)}]{marsden2001discrete}
Marsden, J.E. and West, M. (2001).
\newblock Discrete mechanics and variational integrators.
\newblock \emph{Acta Numerica}, 10, 357--514.

\bibitem[{Modin and S{\"o}derlind(2011)}]{Modin2011}
Modin, K. and S{\"o}derlind, G. (2011).
\newblock Geometric integration of hamiltonian systems perturbed by {R}ayleigh
  damping.
\newblock \emph{BIT Numerical Mathematics}, 51(4), 977--1007.

\bibitem[{Ober-Bl\"obaum(2017)}]{O14}
Ober-Bl\"obaum, S. (2017).
\newblock Galerkin variational integrators and modified symplectic
  {R}unge--{K}utta methods.
\newblock \emph{IMA Journal of Numerical Analysis}, 37(1), 375--406.

\bibitem[{Ober-Bl\"obaum et~al.(2011)Ober-Bl\"obaum, Junge, and
  Marsden}]{ObJuMa10}
Ober-Bl\"obaum, S., Junge, O., and Marsden, J.E. (2011).
\newblock Discrete mechanics and optimal control: an analysis.
\newblock \emph{Control, Optimisation and Calculus of Variations}, 17(2),
  322--352.

\bibitem[{Ober-Bl{\"o}baum and Saake(2015)}]{ober2015construction}
Ober-Bl{\"o}baum, S. and Saake, N. (2015).
\newblock Construction and analysis of higher order galerkin variational
  integrators.
\newblock \emph{Advances in Computational Mathematics}, 41(6), 955--986.

\bibitem[{Patrick and Cuell(2009)}]{patrick2009error}
Patrick, G.W. and Cuell, C. (2009).
\newblock Error analysis of variational integrators of unconstrained lagrangian
  systems.
\newblock \emph{Numerische Mathematik}, 113(2), 243--264.

\bibitem[{Reich(1994)}]{Reic94}
Reich, S. (1994).
\newblock Momentum conserving symplectic integrations.
\newblock \emph{Physica D}, 76(4), 375--383.

\bibitem[{Vermeeren(2017)}]{vermeeren2017modified}
Vermeeren, M. (2017).
\newblock Modified equations for variational integrators.
\newblock \emph{Numerische Mathematik}, 137, 1001--1037.

\bibitem[{Wenger et~al.(2017)Wenger, Ober-Bl{\"o}baum, and
  Leyendecker}]{Wenger2017}
Wenger, T., Ober-Bl{\"o}baum, S., and Leyendecker, S. (2017).
\newblock Construction and analysis of higher order variational integrators for
  dynamical systems with holonomic constraints.
\newblock \emph{Advances in Computational Mathematics}, 43(5), 1163--1195.

\end{thebibliography}

\end{document}